\documentclass[12pt,english]{article}

\usepackage{mathptmx}
\usepackage{geometry}
\usepackage{xcolor}
\usepackage{array}
\usepackage{bbding}
\usepackage{multirow}
\usepackage[fleqn]{amsmath}
\usepackage{amssymb}
\usepackage{amsthm}
\usepackage{graphicx}
\usepackage{natbib}
\usepackage{lscape}
\usepackage[hyphens]{url}
\usepackage[hidelinks]{hyperref}
\usepackage{comment}
\usepackage{bm}
\usepackage[title]{appendix}
\usepackage{longtable}
\usepackage{mathtools}
\usepackage{chngcntr}
\usepackage{authblk}
\usepackage{babel}
\usepackage{dsfont}
\usepackage{subcaption}
\usepackage{fancyhdr}
\usepackage{cleveref}
\usepackage{abstract}

\makeatletter
\allowdisplaybreaks
\date{}

\fancypagestyle{plain}{%
	\fancyhf{} 
	
	\lhead{\footnotesize \ \ \textcolor{gray}{\textbf{}}}
	\chead{\footnotesize \textcolor{gray}{\textit{\shortauthors -- \runningtitle}}}
	\rhead{\footnotesize \textcolor{gray}{\thepage}}
}
\makeatletter
\def\blfootnote{\xdef\@thefnmark{}\@footnotetext}
\makeatother

\makeatletter
\def\titlepageext{
	\begin{center}	
	{\parindent0pt
		\rule{0.9\textwidth}{1pt}
		\begin{minipage}[t]{0.25\textwidth}
			\small {\it Keywords:}\\
			\keyword
		\end{minipage}%
		\hspace{3mm}
		\begin{minipage}[t]{0.6\textwidth}
			\small \abstract
		\end{minipage}%
		
		\rule{0.9\textwidth}{2pt}
	}
	\end{center}

	\blfootnote{* Corresponding author. E-mail address: \href{mailto:\corresemail}{\corresemail}.}
}
\makeatother

\usepackage[mathlines, pagewise]{lineno}
\usepackage{etoolbox} 

\newcommand*\linenomathpatchAMS[1]{%
	\expandafter\pretocmd\csname #1\endcsname {\linenomathAMS}{}{}%
	\expandafter\pretocmd\csname #1*\endcsname{\linenomathAMS}{}{}%
	\expandafter\apptocmd\csname end#1\endcsname {\endlinenomath}{}{}%
	\expandafter\apptocmd\csname end#1*\endcsname{\endlinenomath}{}{}%
}

\expandafter\ifx\linenomath\linenomathWithnumbers
\let\linenomathAMS\linenomathWithnumbers
\patchcmd\linenomathAMS{\advance\postdisplaypenalty\linenopenalty}{}{}{}
\else
\let\linenomathAMS\linenomathNonumbers
\fi

\linenomathpatchAMS{gather}
\linenomathpatchAMS{multline}
\linenomathpatchAMS{align}
\linenomathpatchAMS{alignat}
\linenomathpatchAMS{flalign}

\nolinenumbers

\newtheorem{theorem}{Theorem}

\usepackage{algorithm}
\usepackage[noend]{algpseudocode}

\makeatletter
\def\BState{\State\hskip-\ALG@thistlm}
\makeatother

\usepackage{cleveref}
\usepackage{bbm}
\newcommand{\argmin}{\mathop{\mathrm{arg\,min}}\limits}

\newtheorem{corollary}{Corollary}

\title{Traffic Adaptive Moving-window Service Patrolling for Real-time Incident Management during High-impact Events}

\def\shortauthors{Lei et al.}
\def\runningtitle{Traffic Adaptive Moving-window Patrolling Algorithm}

\author[a]{Haozhe Lei}
\author[a]{Ya-Ting Yang}
\author[a]{Tao Li}
\author[b,$\ast$]{Zilin Bian}
\author[b]{Fan Zuo}
\author[a]{Sundeep Rangan}
\author[b]{Kaan Ozbay}

\affil[a]{Department of Electrical and Computer Engineering, New York University, United States of America}
\affil[b]{Department of Civil and Urban Engineering, New York University, United States of America}

\def\corresemail{zb536@nyu.edu }

\def\abstract{This paper presents the Traffic Adaptive Moving-window Patrolling Algorithm (TAMPA), designed to improve real-time incident management during major events like sports tournaments and concerts. Such events significantly stress transportation networks, requiring efficient and adaptive patrol solutions. TAMPA integrates predictive traffic modeling and real-time complaint estimation, dynamically optimizing patrol deployment. Using dynamic programming, the algorithm continuously adjusts patrol strategies within short planning windows, effectively balancing immediate response and efficient routing. Leveraging the Dvoretzky–Kiefer–Wolfowitz inequality, TAMPA detects significant shifts in complaint patterns, triggering proactive adjustments in patrol routes. Theoretical analyses ensure performance remains closely aligned with optimal solutions. Simulation results from an urban traffic network demonstrate TAMPA's superior performance, showing improvements of approximately 87.5\% over stationary methods and 114.2\% over random strategies. Future work includes enhancing adaptability and incorporating digital twin technology for improved predictive accuracy, particularly relevant for events like the 2026 FIFA World Cup at MetLife Stadium.
}

\def\keyword{High-impact event management, service patrol, dynamic programming, adaptive graph}

\begin{document}
\maketitle
\vspace{-5mm}
\titlepageext
\vspace{-10mm}
\section{Introduction}
\subsection{Motivation}
Organizing high-impact events, such as sports tournaments, festivals, and concerts, presents substantial social, economic, and transportation challenges. These events can place immense pressure on transportation infrastructure, security protocols, and public services, particularly in regions that are already congested and economically vital, such as the New York-New Jersey (NYNJ) or Los Angeles (LA) metropolitan areas. Both regions attract large, diverse crowds as tourists from across state lines and around the world, further complicating special event management logistics. A prime example of this challenge is the hosting of mega-events, such as the FIFA World Cup \cite{ardemagni_security_2022} and the Olympics \citep{tokyogov2022transportation,harrison_securing_2021}. These events are expected to bring a massive influx of visitors, creating major demands on the region's infrastructure and services. 

One of the most pressing challenges during the operations of mega-events such as the World Cup or Olympics is managing the safe and efficient movement of the high density of vehicles and people concentrated in event venues, which creates elevated risks such as violations of security and safety protocols. A central issue is how to allocate limited service patrol resources effectively to respond to the maximum number of complaints rapidly. This challenge is further complicated by the dynamic and uncertain nature of complaints, which can occur unpredictably and vary across the network. Event management must balance addressing immediate complaints, ensuring crowd safety, and maintaining smooth operations, all while dealing with constrained resources and rapidly changing conditions.


To address these needs, service patrols by police and Department of Transportation service trucks must be strategically deployed around event venues to assist tourists and locals, respond to their complaints, and maintain efficient movement of both vehicular and pedestrian traffic flow. During mega-events, the complaints can be highly stochastic, dynamic, and random in both spatial and temporal dimensions, making it extremely challenging for service patrol agents to react and respond rapidly. Moreover, given that service patrol vehicles usually have higher rights of way that are generally available to incident/emergency response vehicles, developing an adaptive strategy to accommodate the constraints of traffic network topology while providing optimized routing plans that leverage their higher rights is critical. Without an efficient strategy, the ability to manage high-density vehicular and pedestrian flow and uphold safety standards is severely compromised.

\begin{figure}[h]
    \centering
    \includegraphics[width=0.8\linewidth]{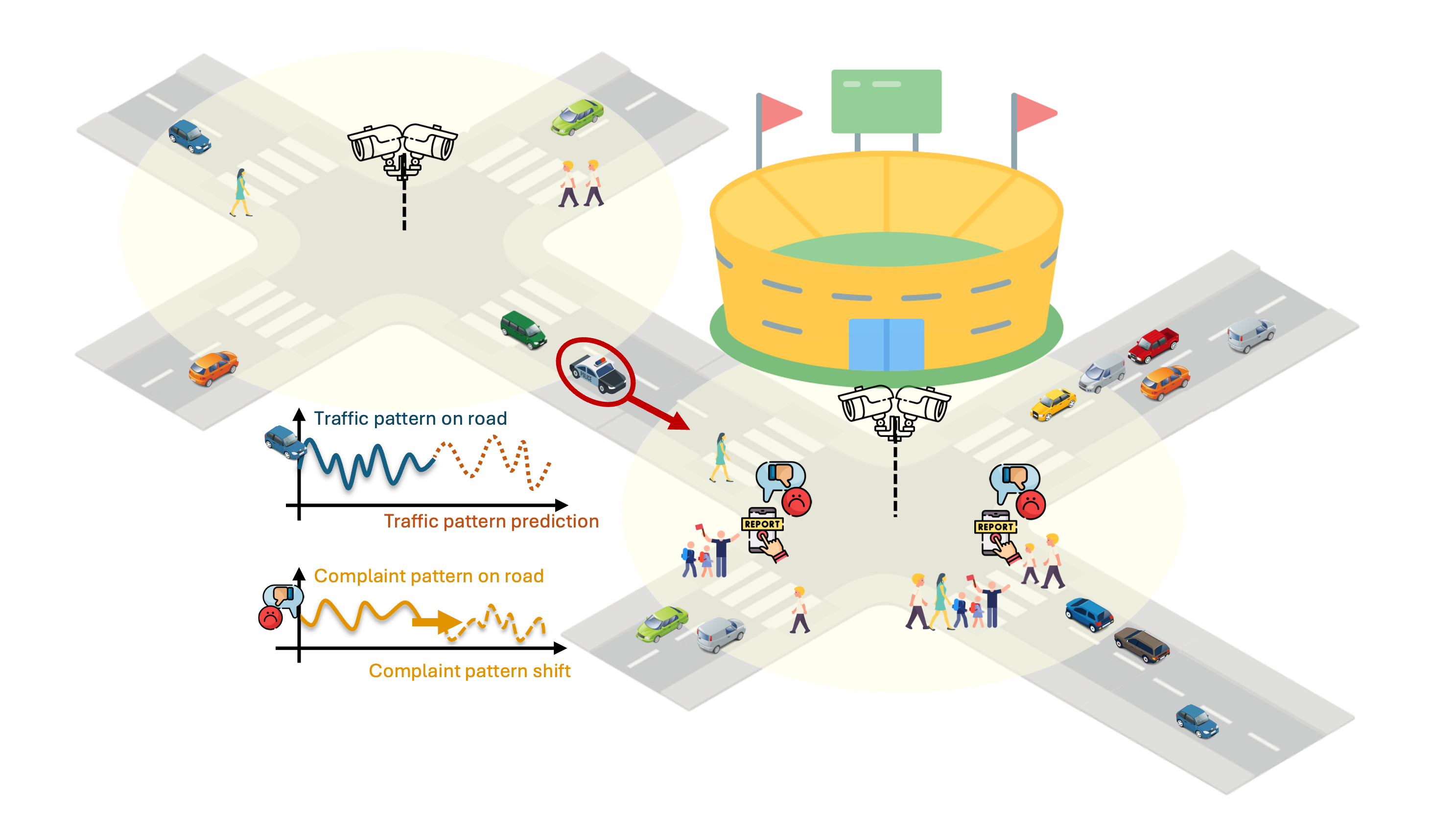}
    \caption{An illustration of TAMPA: the camera captures traffic patterns on the road to assist in generating the prediction of traffic states given incidental situations, while crowd-sourcing collectors (e.g., mobile phone applications) generate user-reported complaint patterns, enabling the detection of shifts in complaints.}
    \label{fig:intro_map}
\end{figure}

In response to the above challenges, we propose a novel traffic patrol algorithm, named \emph{Traffic Adaptive Moving-Window Patrolling Algorithm (TAMPA)}, tailored for faster response to complaints that disrupt traffic while managing high-impact events. The framework addresses non-recurrent congestion resulting from accidents, breakdowns, and other disruptions. Specifically, we incorporate real-time user feedback—treated here as “complaints” or reported incidents—and adapt the patrol route to achieve faster incident clearance. As shown in Fig.\ref{fig:intro_map}, TAMPA leverages existing traffic monitoring and interacting systems (e.g., cameras and mobile phone applications) to sense traffic conditions and complaints in real-time. Subsequently, it optimizes patrol scheduling to minimize travel time to emerging traffic disruptions, delivering swift assistance and improving overall system reliability.

\vspace{-5mm}

\subsection{Our Contributions}
We aim to address two main questions
\begin{itemize}
\item How can service patrol vehicles improve efficiency by balancing the need to address immediate complaints across the transportation network with the constraints imposed by their routing resources and the dynamic nature of network-wide traffic conditions?
\item When facing unexpected complaints (e.g., incidents/accidents), how can the decision maker adapt its service patrol plan in real-time using online observations and change its predetermined patrol route in real-time, leveraging higher rights of way that are generally available to incident /emergency response vehicles (e.g., making U-turns or driving on the shoulder or even against traffic)?
\end{itemize}

To address the first question, the proposed TAMPA integrates a traffic pattern predictor and an estimator for network-wide edge complaint distributions. Using dynamic programming, we determine the optimal strategy for each planning window and continuously calibrate both the predictor and the estimator based on real-time observations of the traffic network during its real-time execution. For the second question, we propose an adapted patrol graph to accommodate higher-priority rights-of-way during incidental and emergency situations. This dynamic approach enables service patrol vehicles to adjust their patrol plans in real-time, bypassing standard traffic rules within the constraints of the traffic network's topological structure. Leveraging the Dvoretzky–Kiefer–Wolfowitz testing \citep{DKW}, TAMPA triggers immediate adaptation of the patrol graph and initiates a new planning window whenever the observed distribution of complaints significantly deviates from its prior empirical distribution. To justify the event-triggered graph update in TAMPA, we prove in \Cref{th:bounded_cost} that the patrol performance increment due to the new planning window is upper-bounded by the total variation distance between the prior and observed distributions multiplied by the planning horizon. If the deviation is below a given threshold, the previous patrolling plan still delivers a suboptimal performance. Hence, TAMPA does not update the current patrolling plan to reduce computation overhead, achieving online adaptive and efficient patrol planning. Utilizing similar proof techniques, we further prove in \Cref{th:bounded_cost_true} that the suboptimality gap of the proposed TAMPA with respect to the optimal patrolling plan is upper-bounded by the total variation between the observed and actual distribution multiplied by the planning horizon. Finally, we validate our approach through simulation experiments based on real-world scenarios, demonstrating its effectiveness in environments with non-stationary traffic patterns and dynamic complaint distributions (around $87.5\%$ improvement compared with stationary strategy and $114.2\%$ improvement compared with random patrol strategy).
\vspace{-5mm}

\section{Related Works}
Non-recurrent congestion due to incidents, which is responsible for more than 50 percent of congestion, can be reduced using efficient incident management strategies such as the use of service patrols for rapid and timely detection and response (\cite{ozbay1999incident,ma2009harnessing,baykal2009modeling}). Most previous patrol studies focus on patrol for crime rate reduction and crime incident response time minimization. In this context, \cite{samanta2022literature} has categorized such police patrolling problems into three categories:  \textit{District design} \cite{liberatore2020police} that partitions the whole area into several regions termed as districts in order to help balance the workload depending on the crime density in each district; \textit{Resource allocation} that determines how to allocate and schedule the available resources like service patrol vehicles or manpower to the predefined hotspots \cite{keskin2012analysis} or throughout the entire network to maximize coverage \cite{oghovese2014optimal} or improve agent visibility \cite{adler2014location}, and \textit{Route design} \cite{azimi2016modeling,aguirre2012evolutionary} that designs hard-to-predict patrol paths and coordinates patrollers to cover the maximum area. 

In addition to crime-reduction objectives, service patrol can also benefit other applications. As for speed and red light violation prevention, \cite{dabaghchian2019intelligent} proposes a learning algorithm that assigns two police officers to different intersections to achieve optimal traffic violation prevention at the end of the time horizon. \cite{adler2014location} utilizes a series of linear programming problems for the traffic police routine patrol vehicle assignment that intends to create a presence that acts as a deterrence at intersections and issuing tickets to traffic regulation offenders. In applications for traffic accident management \cite{yazici2015evaluation} or emergency evacuation \cite{li2015evacuation}, \cite{ozbay2004probabilistic} utilizes mathematical programming models with probabilistic constraints to address the dispatching problem in incident response as well as allocate service vehicles to each depot. \cite{wang2021efficient} considers both online and offline models for more adaptive patrol. The recent study \cite{repasky2024multi} develops a reinforcement learning algorithm that takes both preventive regular service patrol and being dispatched to serve emergency incidents, while \cite{wang2024improving} proposes a simulation-based framework to help plan the optimum patrol routes based on available safety service patrol resources and predicted incidents.

However, when it comes to patrol for (mega) event management, there are very limited studies focused on this specific application. For example, \cite{camacho2015decision} proposes a patrol decision support system aimed at crime reduction but only briefly suggests that their model could be extended in future work to include variables related to special events (e.g., demonstrations, sports events, and parades).  This work, therefore, seeks to fill this gap by exploring the role of service patrols in supporting and enhancing (mega) event management.

\section{Problem Formulation}
Denote the patroller is responsible for designing the patrol route plan for the service patrol vehicle and defining the patrol procedure over a finite timeline, typically within a single day (an appropriate duration for any service patrol vehicle). This timeline is represented as a discrete set $t \in T := \{0, 1, \ldots, |T|\}$. We assume that the traffic conditions evolve over time $t$, meaning that all components listed below are dependent on $t$. Let the patrol graph $\mathcal{G}_t = \langle \mathcal{V}_t, \mathcal{E}_t \rangle$ belonging to the set of graphs $\mathcal{Q}$ represent patrol-designated areas associated with road connections in a traffic network at time $t$, where $\mathcal{V}_t$ is the set of designated areas (e.g., intersections, bus stations, traffic monitoring cameras), and $\mathcal{E}_t \subseteq \mathcal{V}_t \times \mathcal{V}_t$ is the set of directed edges representing road connections from node $u$ to node $v$. We define the set of edge lengths as $
\mathcal{L}_t = \{ l^{(i,j)} \mid (i,j) \in \mathcal{E}_t \},$ where $ l^{(i,j)}\in \mathbb{R}_{> 0} $ is the length of edge $(i,j)$ at time $t$ satisfying $l(i,j) = l(j,i), \text{ for all } (i,j) \in \mathcal{E}_t.$ For the use in the following sections, we also define the neighbor of a node $v$, including itself, as $N(v) = \{u \in \mathcal{V} \mid (v, u) \in \mathcal{E}\} \cup \{v\}.$ And the spatial coordinates of each node as $\{x_v\}_{v\in\mathcal{V}_t}$.

Note that although $\mathcal{G}_t$ is a directed graph, the length of edge $(i,j)$ is equal to the length of edge $(j,i)$, ensuring consistency in edge lengths irrespective of their direction. Regularly, $\mathcal{G}_{t}=\mathcal{G}_{t+1}$, we will talk about when it will change in \Cref{sec:GAP}. Denote the sequence of travel times as $\{\{\mu^e_t\}_{e \in \mathcal{E}_t}\}_{t \in T}$, where each $\mu^e_t \in \mathcal{U}\subset \mathbb{R}$ represents the travel time on edge $e$ at time $t$. Let the minimum travel time $\{mtt^e\}_{e \in \mathcal{E}_t}$ be determined by road capacity and speed limits; then, we have the constraint $\mu^e_t \geq mtt^e, \, \forall e \in \mathcal{E}_t, \, \forall t \in T$. Suppose there is a sequence of complaints, $\{\{c^e_t\}_{e \in \mathcal{E}_t}\}_{t \in T}$, where the number of complaints $ c^e_t $ is a random variable. We assume that $ c^e_t $ is distributed according to the distribution $ \mathcal{C}_t^e $, which is supported on the set of non-negative real numbers, i.e.,$c^e_t \sim \mathcal{C}_t^e$, where $ \mathcal{C}_t^e$ is supported on $\mathbb{R}_{\geq 0}$ and varies over time and location.

\begin{figure}[H]
    \centering
    \includegraphics[width=1.0\linewidth]{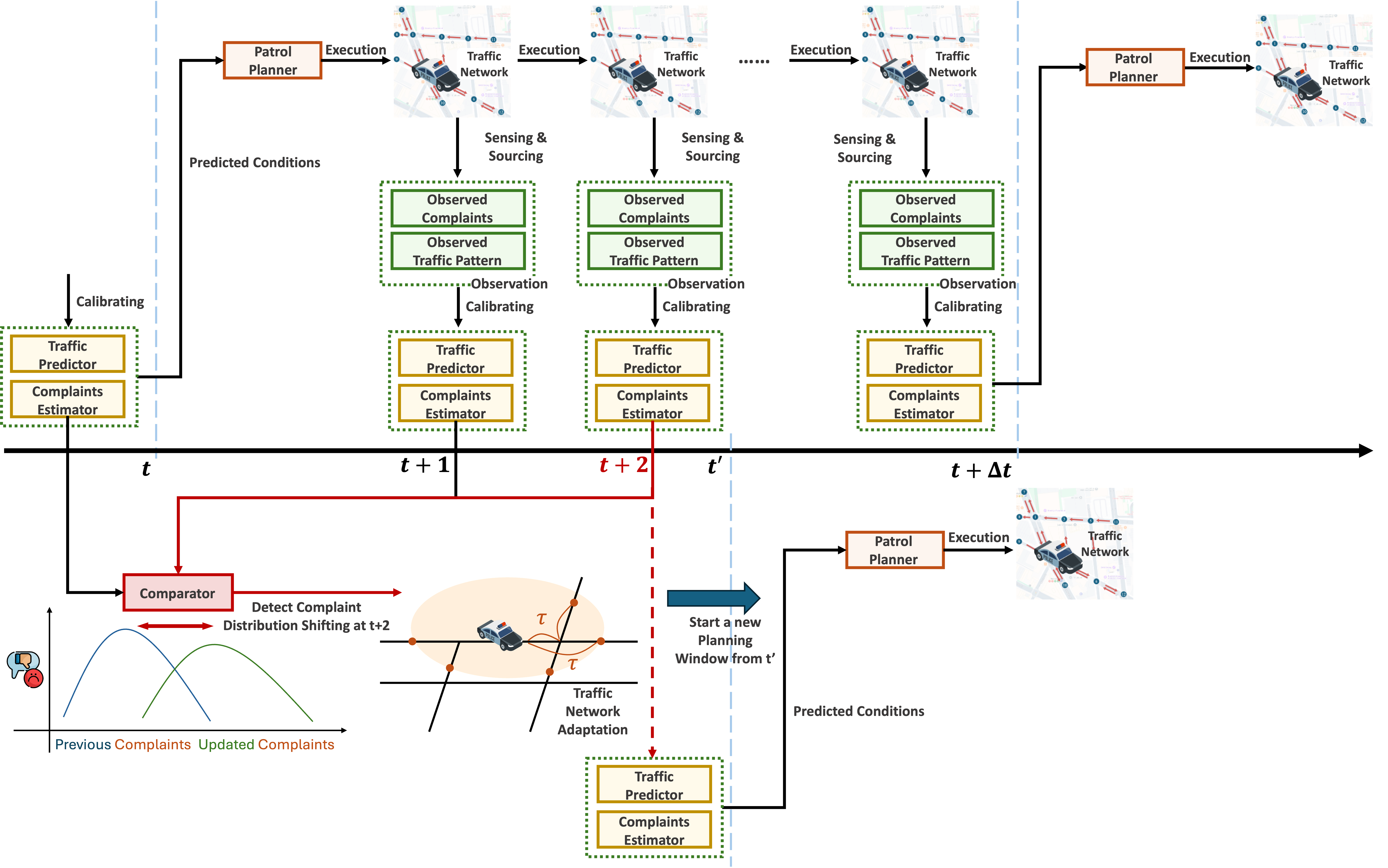}
    \caption{The illustration of TAMPA that combines online solution of the planning window MDP problem (above) in \Cref{sec:pwmdp} and the global adaptation process (below) in \Cref{sec:GAP}. The TAMPA provides the patrol command to the service patrol vehicle at each time step. In solving the planning window MDP problem, the traffic predictor and the complaints estimator provide predicted conditions for the service patrol planner. The planner obtains an optimal action, and the service patroller executes it in the traffic network. The sensing via camera to observe traffic conditions and the crowd-sourcing collector to generate user-reported complaints in the network, then update the traffic predictor and the complaints estimator at each time step until the start of the next planning window. The global adaptation process helps the service patrol vehicle proactively adapt to the variation of the complaints conditions. The comparator will decide if complaints distribution shifting happens at each time step. If shifting happens, the patrol graph will be adapted, and a new planning window will be initiated immediately.}
    \label{fig:algo_graph}
    \vspace{-5mm}
\end{figure}

The service patrol vehicle operates under two distinct statuses: \textbf{inspecting} and \textbf{commuting}. The \textbf{inspecting} status indicates that the patroller remains stationed at a node, whereas the \textbf{commuting} status entails traveling between different nodes. We assume that the former status allows the patrol vehicle to satisfy complaints on its neighboring edges, restricted by the patroller's inspecting threshold $\zeta \in \mathbb{R}_{>0}$. This limitation reflects the practical constraint that the patroller can inspect only a finite distance, thereby addressing complaints solely on edges where the edge length does not exceed the allowable inspecting threshold. The latter status enables the patroller to satisfy complaints on the edge it is traveling on, representing the patroller's ability to address complaints during its traversal. 

Define a set of time points $\{\eta_0, \eta_1, \dots, \eta_{\text{end}}\} \subseteq T,$
where each $\eta_g$ represents an instant at which the patroller is located at a node (as opposed to traveling along an edge). Let $v_{\eta_g}$ denote the patroller's node at time $\eta_g$, and let $v_{\eta_{g+1}}$ denote the patroller's subsequent destination node, we obtain the trajectory of the patroller is $\mathcal{T}:=\{v_{\eta_0}, v_{\eta_1}, \ldots, v_{\eta_{\text{end}}}\}$.

Assuming the patroller has perfect information for the entire day, any stochastic distribution measurement can be treated as a known (deterministic) realization at each time point. Suppose the patroller's utility at time $\eta_g$ is given by $ r_g = r_g\bigl(v_{\eta_g}, v_{\eta_{g+1}}\bigr)$. We introduce a weight parameter $\lambda \ge 0$ and define the patroller's utility in two cases:

\begin{itemize}
    \item \textbf{Inspecting (stationary):} When $v_{\eta_g} = v_{\eta_{g+1}}$,
    \[
        r_g 
        \;=\; (1-\lambda)\sum_{u \in N\bigl(v_{\eta_g}\bigr)} 
        c_{\eta_g}^{\bigl(v_{\eta_g}, u\bigr)} 
        \,\cdot\, \min\!\Bigl(1,\; \frac{\zeta}{l^{(v_{\eta_g},u)}}\Bigr).
    \]
    \item \textbf{Commuting (moving):} When $v_{\eta_g} \neq v_{\eta_{g+1}}$,
    \[
        r_g 
        \;=\; -\,\lambda\,\mu_{\eta_g}^{\bigl(v_{\eta_g}, v_{\eta_{g+1}}\bigr)}
        \;+\;(1-\lambda)\,c_{\eta_g}^{\bigl(v_{\eta_g}, v_{\eta_{g+1}}\bigr)}.
    \]
\end{itemize}

\noindent
Our goal is to determine an optimal patrolling trajectory that maximizes the total utility over the sequence of time points:
\begin{align}
  Q^* = \max_{\mathcal{T}}
  \sum_{g=0}^{\text{end}-1} 
  r_g\bigl(v_{\eta_g}, v_{\eta_{g+1}}\bigr).
  \label{eq:global_opt}
\end{align}

\section{Traffic Adaptive Moving-window Patrolling Algorithm (TAMPA)}

In practice, the future cannot be known with certainty, and thus, the realization of a probability distribution cannot be determined in advance. Therefore, in this section, we propose a method called TAMPA shown in \Cref{fig:algo_graph} to address the patrolling problem under these realistic constraints.

\subsection{Planning Window MDP Problem}\label{sec:pwmdp}
To simplify the model in deciding the optimal strategy for the patroller, we define a constant fine-tuned time interval $\tau \in \mathbb{R}_{>0}$ as the length of each decision time slot $k \in K := \{0, 1, \cdots, |K|\}$. The planning window $W_{\mathbf{t}, \tau, |K|}$ is defined over the finite timeline $T$ and initiates at the start time point $ \mathbf{t} \in T$. It comprises $ |K| $ discrete decision-making points, each separated by the fixed time interval $ \tau $, where $ \tau $ is an integer multiple of the base time unit of $ T $. The total duration, or horizon length, of the planning window is given by $ H = |K| \cdot \tau $, encompassing the time span from $ \mathbf{t} $ to $ \mathbf{t} + H $. Formally, the planning window can be expressed as $W_{\mathbf{t}, \tau, |K|} = \left\{ \mathbf{t} + i \cdot \tau \;\middle|\; i = 0, 1, 2, \ldots, |K|-1 \right\},$ which delineates a structured sequence of decision points that facilitate dynamic and finite-horizon patrol planning by allowing periodic updates based on evolving traffic conditions within the specified temporal scope.

Since traffic conditions vary over time and cannot be completely predicted in a protracted time span, we consider a predictor that is sufficiently accurate within the limited time horizon $H$. For a given prediction starting time $\mathbf{t}$, denote the travel time predictor as $\Xi(W_{\mathbf{t}, \tau, |K|}):=\{\{\hat{\mu}^e_k\}_{e \in \mathcal{E}}\}_{k \in \{t, t+\tau,\cdots,t+|K|\tau \}}.$

For simplicity, we consider the approximation of complaint distributions
$\{\hat{c}_k^e \sim \hat{\mathcal{C}}_{\mathbf{t}}^{e}\}_{e \in \mathcal{E}}$, 
where each $\hat{\mathcal{C}}_{\mathbf{t}}^{e}$ is a discrete probability distribution 
over the set of non-negative integers. Concretely, for each edge 
$e \in \mathcal{E}$, we have
\begin{align*}
\hat{\mathcal{C}}_{\mathbf{t}}^{e}(n) = \mathbb{P}\bigl(\hat{c}^e = n\bigr), 
\quad \text{for } n \in \{0,1,2,\dots\},
\end{align*}
with $\sum_{n=0}^{\infty}\hat{\mathcal{C}}_{\mathbf{t}}^{e}(n) = 1.$
A slot complaint for a given estimation start time $\mathbf{t}$ 
is defined as $\hat{c}_k^e = \sum_{t = \mathbf{t}}^{\mathbf{t} + k} \hat{c}^e_t.$ We then define the estimation of a sequence of slot complaints as
\begin{align*}
\Phi\bigl(W_{\mathbf{t}, \tau, |K|}, \{\hat{\mathcal{C}}_{\mathbf{t}}^{e}\}_{e \in \mathcal{E}}\bigr)
:= \bigl\{\{\hat{c}^e_k\}_{e \in \mathcal{E}}\bigr\}_{k \in K}.
\end{align*}

For any planning window $W_{\mathbf{t}, \tau, |K|}$, with given $\{\hat{\mathcal{C}}_{\mathbf{t}}^{e}\}_{e \in \mathcal{E}_{\mathbf{t}}}$,  $\Xi(\cdot)$, $\Phi(\cdot)$, and an initial condition $\rho = \{v_0 = v_{\mathbf{t}}, \mathcal{V} = \mathcal{V}_{\mathbf{t}}, \mathcal{E}= \mathcal{E}_{\mathbf{t}}, \mathcal{L}=\mathcal{L}_{\mathbf{t}}\}$ (assume the patrol graph remains stationary within any moving window), we denote a Markov Decision Process (MDP) for this finite discrete step patrol planning problem as $\langle \mathcal{S}, \mathcal{A}, \mathbf{P}, \mathbf{R}, \rho \rangle$. Let $s_k \in \mathcal{S}$ be the state, and $s_k = \{v_k\}$, so that $\mathcal{S} = \mathcal{V}$ is the node that the patroller is visiting. To reduce the problem's complexity and align with the complaints estimation method mentioned in the previous section, without loss of generality, we assume that the time required for implementing each action is also the given constant time $\tau$ introduced in the previous paragraph. 
Let $\mathcal{P}_{ij}$ denote the set of all feasible paths from node $i$ to node $j$ in the graph $\mathcal{G} = \langle \mathcal{V}, \mathcal{E}\rangle$, where a path $p \in \mathcal{P}_{ij}$ is a sequence of edges connecting $i$ to $j$ without revisiting any node. In a weighted graph with travel time $\hat{\mu}^e_k$ on each edge $e \in \mathcal{E}$, we define the shortest path $p^*_k(i,j)$ and its distance function $d^*_k(i,j)$ under the travel time $\hat{\mu}^e_k$ as 
\begin{align*}
p^*_k(i,j)=\argmin_{p \in \mathcal{P}_{ij}} \;\sum_{\,e \in p}\; \hat{\mu}^e_k, \quad d^*_k(i,j) = \sum_{e \in p^*_k(i,j)} \hat{\mu}^e_k,
\end{align*}
In particular, if $i=j$ we define $d^*_k(i,j) = 0$.

Then, the action for every decision-making step is $a_k \in \mathcal{A}_k(s_k)$ such that $\mathcal{A}_k(s_k) = \{a \in \mathcal{V} \mid d^*_k(s_k,a)\leq \tau\}$ where the transition function $\mathbf{P}$ now becomes a deterministic indicator such that $s_{k+1} = a_k$ with a probability equal to $1$.

Our objective is to find an optimal control sequence that minimizes the cumulative travel time and maximizes the satisfaction of cumulative complaints. In this case, we consider a routing cost
\begin{align}
    r_{tt}(s_k, a_k) = 
    \begin{cases}
        0, & \text{if } s_k = a_k \\
        d^*_k(s_k,a_k), & \text{otherwise}
    \end{cases},
\label{eq:travel_time_cost}
\end{align}
And a complaint satisfaction cost as
\begin{align}
    r_c(s_k, a_k) = 
    \begin{cases}
        \sum_{u \in N(s_k)} \hat{c}_k^{(s_k, u)} \cdot \min\left(1, \frac{\zeta}{l^{(s_k, u)}}\right), & \text{if } s_k = a_k \\
        \sum_{e\in p^*_k(s_k,a_k)}\hat{c_k}^{e}, & \text{otherwise}.
    \end{cases},
\label{eq:complaint_cost}
\end{align}
We assume the reward function $\mathbf{R}$ is equal to the weighted sum of negative routing cost and complaints satisfactory such that:
\begin{align}
    \mathbf{R}(s_k,a_k) = -\lambda\cdot r_{tt}(s_k,a_k)+ (1-\lambda)\cdot r_c(s_k, a_k),
\label{eq:reward_function}
\end{align}
where $\lambda\geq 0$ is the reward weight. Let $s_0=\{v_0\}$ be the initial state where the service patrol agent starts. Its objective function is to find an action sequence $\hat{\mathbf{a}}^*=\{a_0,a_1,\cdots,a_{|K|-1}\}$ that can maximize the cost function that is:
\begin{align}
    J_{\mathbf{t}}(s_{\mathbf{t}})=\mathbb{E}\left\{\sum_{k=0}^{|K|-1}\mathbf{R}(s_{\mathbf{t}+k},a_{\mathbf{t}+k})\right\}.
\label{eq:objective}
\end{align}

Then $a_0\in\hat{\mathbf{a}}^*$ is the optimal action $a^*(W_{\mathbf{t}, \tau, |K|})$ for this planning window. Consider the optimal action $a^*(W_{\mathbf{t}, \tau, |K|})$ obtained by solving the patrolling problem within the planning window. Executing this action incurs a cost of $\Delta_{\mathbf{t}} = \mu^{(v_{\mathbf{t}}, a^*(W_{\mathbf{t}, \tau, |K|}))}_{\mathbf{t}}$. Consequently, the subsequent planning window should commence at time $\mathbf{t}' = \mathbf{t} + \Delta_{\mathbf{t}}$.

\begin{algorithm}[H]
\caption{Finite-Horizon Patrolling Decision Procedure}
\label{alg:patrol-decision}
\begin{algorithmic}[1]
\Require $W_{\mathbf{t}, \tau, |K|}$, $\{\hat{\mathcal{C}}_{e}\}_{e}$, $\Xi(W_{\mathbf{t}, \tau, |K|})$, $\Phi(W_{\mathbf{t}, \tau, |K|}, \{\hat{\mathcal{C}}_{e}\})$, initial $\rho=\{v_0=v_{\mathbf{t}}, \mathcal{V}, \mathcal{E}, \mathcal{L}\}$, and $\lambda$
\State Define MDP $\langle \mathcal{S}, \mathcal{A}, \mathbf{P}, \mathbf{R}, \rho \rangle$ with $\mathcal{S}=\mathcal{V}$, $\mathcal{A}(s_k)$ per definition in text.
\State Compute costs $r_{tt}(s_k,a_k)$ \eqref{eq:travel_time_cost} and $r_c(s_k,a_k)$ \eqref{eq:complaint_cost}.
\State Determine reward function $\mathbf{R}(s_k,a_k)$ \eqref{eq:reward_function}.
\State Solve finite-horizon MDP to find $\hat{\mathbf{a}}^*$ that maximizes $J_{\mathbf{t}}(s_{\mathbf{t}})$ \eqref{eq:objective}.
\State Let $a^*(W_{\mathbf{t}, \tau, |K|}) = a_0$ from $\hat{\mathbf{a}}^*$.
\State Compute $\Delta_{\mathbf{t}}=\mu^{(v_{\mathbf{t}}, a^*(W_{\mathbf{t}, \tau, |K|}))}_{\mathbf{t}}$ and set next planning start $\mathbf{t}' = \mathbf{t}+\Delta_{\mathbf{t}}$.
\end{algorithmic}
\end{algorithm}

\subsection{Global Adaptation Process}\label{sec:GAP}
Since stationary predictions cannot be guaranteed to be 100\% accurate, that is, $\Xi\left(W_{\mathbf{t}, \tau, |K|}\right) \neq \left\{\{\mu^e_k\}_{e \in \mathcal{E}}\right\}_{k \in \{t, t+\tau, \ldots, t + |K|\tau\}}$ and $\Phi\left(W_{\mathbf{t}, \tau, |K|}, \{\hat{\mathcal{C}}_{e}\}_{e \in \mathcal{E}}\right) \neq \left\{\{c^e_k\}_{e \in \mathcal{E}}\right\}_{k \in K},$ the traffic network topology $\mathcal{G}$ may undergo dynamic changes. For instance, such uncertainties might necessitate the patroller to perform a U-turn to reschedule the next destination; thus, a new turning node needs to be added, leading to the updating of $\mathcal{G}$. Consequently, there is a need to model an online adaptation procedure for the patrol problem.

Following our previous work \citep{ttcx}, we assume the traffic condition 
$o_t:= \left\{\{\mu^e_{t}\}_{e \in \mathcal{E}_t}, \{c^e_{t}\}_{e \in \mathcal{E}_t}\right\}$ is completely observable by traffic management device such as tilting cameras.

Given that the exact distribution $ \mathcal{C}_t^e $ governing the traffic condition $ c^e_t $ on edge $ e \in \mathcal{E}_t $ is unknown, and suppose we observe the sample of edge complaints $c_{t}^{e} \sim \mathcal{C}_{t}^{e}$ at each time $t$. We maintain an empirical approximation $\hat{\mathcal{C}}_{t}^{e}$ of $\mathcal{C}_{t}^{e}$ using the observed samples. Following the Glivenko–Cantelli theorem, let $\hat{\mathcal{C}}_{0}^{e}(n) = \mathbb{P}\bigl(\hat{c}^{e}=n\bigr),  n \in \{0,1,2,\dots\},$
be a discrete prior distribution (with possibly infinite support) that represents our initial belief about $\hat{c}^{e}$.  Suppose this prior corresponds to $M$ ``historical samples.''  At each time step $t\in\{1,2,\ldots\}$, we observe a single data point $    c_{t}^{e}$ (also taking values in $\{0,1,2,\dots\}$).
We wish to update our approximation of the underlying distribution from $\hat{\mathcal{C}}_{t-1}^{e}$ to $\hat{\mathcal{C}}_{t}^{e}$ in an online method. We use $\delta_{c_{t}^{e}}(n)$ to denote the Dirac mass at $c_{t}^{e}$, i.e.,
\begin{align*}
    \delta_{c_{t}^{e}}(n) = 
    \begin{cases}
      1, & \text{if } n = c_{t}^{e}, \\
      0, & \text{otherwise.}
    \end{cases}
\end{align*}

By assigning weight $M$ to the prior distribution and $1$ to each new observation, we define the updated distribution at time $t$ via
\begin{align}
\hat{\mathcal{C}}_{t}^{e}(n)
    \;=\;
    \frac{M + t - 1}{M + t} \;\hat{\mathcal{C}}_{t-1}^{e}(n)
    \;+\;
    \frac{1}{M + t}\;\delta_{c_{t}^{e}}(n).
    \label{eq:update}
\end{align}
Unrolling this recursion from $t=1$ up to $t=|T|$, we obtain the closed-form expression
\begin{align*}
\hat{\mathcal{C}}_{|T|}^{e}(n)
    \;=\;
    \frac{S}{S + |T|}\;\hat{\mathcal{C}}_{0}^{e}(n)
    \;+\;
    \frac{1}{S + |T|}\;\sum_{t\in T} \delta_{c_{t}^{e}}(n).
\end{align*}

According to the Glivenko–Cantelli theorem \cite{tucker1959generalization}, for a large enough $T$,
\begin{align*}
    \sup_{n \in \{0,1,2,\cdots\}}|\hat{\mathcal{C}}_{|T|}^{e}(n)-\mathcal{C}^{e}|\rightarrow 0 \text{ almost surely.}
\end{align*}
where with a slight abuse of symbol, $\mathcal{C}^{e}$ represents a distribution for the complete timeline $T$.

For moderate $T$, the prior $\hat{\mathcal{C}}_{0}^{e}$ can still exert significant influence, effectively pulling the estimate toward one’s initial belief. All in all, since the variation of traffic complaints traditionally cannot be capricious, the approximated distribution $\hat{\mathcal{C}}_t^e$ is approaching the exact distribution $\mathcal{C}_t^e$ via the updating of itself.

To decide whether the updated distribution function $\hat{\mathcal{C}}_{t}^{e}$ of complaints is different from a known prior distribution $\overline{\mathcal{C}}^e$, (let $\overline{\mathcal{C}}^e=\hat{\mathcal{C}}_{0}^{e}$ as the initial setting), we declare that the complaint distribution is shifted happened on edge $e$ if event $E_t^e:=\sup_{n \in \{0,1,2,\cdots\}} \left| \overline{\mathcal{C}}^e(n) - \hat{\mathcal{C}}_{t}^{e}(n) \right| \geq q$ is true. For example, let $q=\sqrt{3/2t}$ as our selected shifting threshold. Refer to the Dvoretzky–Kiefer–Wolfowitz inequality:
\begin{align*}
\mathbb{P} \left\{ \sup_{n \in \{0,1,2,\cdots\}} \left| \overline{\mathcal{C}}^e(n) - \hat{\mathcal{C}}_{t}^{e}(n) \right| \geq q \right\} \leq 2e^{-2tq^2}, \quad \forall s \geq 0.
\end{align*}
The event that two distributions still remain the same happens at most $10\%$. We call it the DKW testing.

Define the network DKW divergence measures $D_t^q$ in \eqref{eq:DKW_index} as a single metric that aggregates the "shift or no-shift" outcomes for all edges in the network into one comprehensive measure. 

\begin{align}
    D_t^q = \prod_{e \in \mathcal{E}_t}\mathbbm{1}\left(E_t^e\right).\label{eq:DKW_index}
\end{align}

Specifically, for each edge $e \in \mathcal{E}_t$, it uses an indicator function $ \mathbbm{1}(E_t^e)$ to determine whether the complaint distribution has shifted beyond the threshold $q$. By taking the product of these indicators across all edges, $D_t^q$ quantifies the overall likelihood that the network’s complaint distributions have significantly diverged from their prior distributions.

The update of $ \overline{\mathcal{C}}^e $ occurs under two distinct scenarios. The first scenario arises when a new planning window commences, at which point we set $ \overline{\mathcal{C}}^e = \hat{\mathcal{C}}_{\mathbf{t}}^{e} $ where $\mathbf{t}$ is the last moving window's start time point. The second scenario occurs if the divergence measures $ D_t^q $ satisfies $ D_t^q \neq 0 $, indicating a significant shift in the complaint distributions. In this case, we update the prior distribution to $ \overline{\mathcal{C}}^e = \hat{\mathcal{C}}_{\mathbf{t}}^{e} $. Consequently, a new planning window is constructed, and a new optimal action is determined immediately, effectively enabling proactive planning ahead of the next planning window.

Following similar proof techniques in  \cite{li2024metastackelberggamerobust}, we prove \Cref{th:bounded_cost} to demonstrate that, with an appropriate choice of $q$, if the network-wide edge distribution satisfies $D^q_t > 0$, then despite updates and differences in the distributions compared to earlier time points, the performance gap in the cost of the last planning window (see \Cref{sec:pwmdp})—using the updated complaints distribution—remains bounded.

Refer to \eqref{eq:objective}, a bit of abuse of notations, denote $J_{\mathbf{t}}[\mathcal{C}](v_{\mathbf{t}})$ as the cost function of the planning window at $\mathbf{t}$ with network-wide edge distributions $\mathcal{C}$. Let $\|P_1 - P_2\|_{\mathrm{TV}}$ be the total variation distance between two distributions $P_1$ and $P_2$.

\begin{theorem}

If $\mathbf{t}$ denotes the start time of the last planning window, and $\mathcal{O}(\cdot)$ denotes the Big-O notation,
\begin{align*}
    \bigl|J_{\mathbf{t}}[\overline{\mathcal{C}}](v_{\mathbf{t}}) \;-\; J_{\mathbf{t}}[\hat{\mathcal{C}}_{t}](v_{\mathbf{t}})\bigr|
    \;\leq\;
    \mathcal{O}\!\Bigl(|K|\sum_{e \in \mathcal{E}_t} \bigl\| \overline{\mathcal{C}}^e - \hat{\mathcal{C}}_{t}^e \bigr\|_{\mathrm{TV}}\Bigr).
\end{align*}
\label{th:bounded_cost}
\end{theorem}

Furthermore, we prove \Cref{th:bounded_cost_true} to show the performance gap is also bounded using this empirical distribution versus using the true distributions 

\begin{corollary}
If $\mathbf{t}$ denotes the start time of the last planning window and $\mathcal{C}_t$ is the true network-wide edge distributions at time $t$
\begin{align*}
    \bigl|J_{\mathbf{t}}[\mathcal{C}_t](v_{\mathbf{t}}) \;-\; J_{\mathbf{t}}[\hat{\mathcal{C}}_{t}](v_{\mathbf{t}})\bigr|
    \;\leq\;
    \mathcal{O}\!\Bigl(|K|\sum_{e \in \mathcal{E}_t} \bigl\| \mathcal{C}^e_t - \hat{\mathcal{C}}_{t}^e \bigr\|_{\mathrm{TV}}\Bigr).
\end{align*}
\label{th:bounded_cost_true}
\end{corollary}
The complete proofs of \Cref{th:bounded_cost} and \Cref{th:bounded_cost_true} are provided in the main paper, as they are omitted here due to space constraints in this extended abstract.

There are two distinct statuses for the patroller when the distribution shifting happens; accordingly, there are two different resolving procedures.

\subsubsection{Inspecting:} No changes are required. The patroller simply initiates the solving of an adaptive patrolling problem within the planning window, starting from the current time $ \mathbf{t}_{\text{new}} = t+1 $. This is represented as $\left\langle W_{\mathbf{t}_{\text{new}}, \tau, |K|}, \{\hat{\mathcal{C}}_{e}\}_{e \in \mathcal{E}_{\mathbf{t}_{\text{new}}}}, \Xi(\cdot), \Phi(\cdot), \rho = \{v_0 = v_{\mathbf{t}_{\text{new}}}, \mathcal{V} = \mathcal{V}_{\mathbf{t}_{\text{new}}}, \mathcal{E}= \mathcal{E}_{\mathbf{t}_{\text{new}}}, \mathcal{L}=\mathcal{L}_{\mathbf{t}_{\text{new}}}\} \right\rangle.$

\subsubsection{Commuting:}\label{sec:commuting_update} In this scenario, when $ t = \mathbf{t}' $, it signifies the beginning of the next planning window; thus, no changes are required. If the shifting occurs during $ t \in (\mathbf{t}, \mathbf{t}') $, it necessitates rescheduling the new planning window on an edge rather than on a node. Let $v_o$ denote the origin node and $v_d$ the destination node. Together, $(v_o, v_d)$ form an O-D pair, which represents the patroller's current commuting route. To update the traffic network graph $\mathcal{G}_t = \langle \mathcal{V}_t, \mathcal{E}_t \rangle$, new nodes are added to split the commuting route $(v_o=v_{\mathbf{t}}, v_d=a^*(W_{\mathbf{t}, \tau, |K|}))$ and its neighbors edges based on the ratio of the travel time. Let the patroller start at node $v_o$ at time $\mathbf{t}$ and move toward the destination node $v_d$, with the estimated time of arrival $\mathbf{t}'$. At the current time $t$, the ratio of time spent to the total travel time is given by 

\begin{align*}
\gamma_t = \frac{t - \mathbf{t}}{\mathbf{t}' - \mathbf{t}}.
\end{align*}
One edge $(v_o,v_d)$, the new node $v_{\text{new}}^{(v_o, v_d)}$ is defined as a point that splits the edge $(v_o, v_d)$. Let the spatial coordinates of $v_o$ and $v_d$ be denoted by $\mathbf{x}_t$ and $\mathbf{x}_d$, respectively. The spatial position of $v_{\text{new}}^{(v_o, v_d)}$ is calculated as:
\begin{align}
\mathbf{x}_{\text{new}}^{(v_o, v_d)} = (1 - \gamma_t)\mathbf{x}_t + \gamma_t\mathbf{x}_d.
\label{eq:posnew}
\end{align}

The updated set of nodes becomes:
\begin{align}
\mathcal{V}^{'}_{t} = \mathcal{V}_t \cup \{v_{\text{new}}^{(v_o, v_d)}\}.
\label{eq:Vupdate}
\end{align}

The edge $(v_o, v_d)$ is then replaced by two new directed edges: $(v_o, v_{\text{new}}^{(v_o, v_d)})$ and $(v_{\text{new}}^{(v_o, v_d)}, v_d)$. The updated set of edges is:
\begin{align}
\mathcal{E}^{'}_{t} = \{\mathcal{E}_t \setminus \{(v_o, v_d)\}\} \cup \{(v_o, v_{\text{new}}^{(v_o, v_d)}), (v_{\text{new}}^{(v_o, v_d)}, v_d)\}.
\label{eq:Eupdate}
\end{align}

Accordingly, the corresponding traffic flow distribution $\{\{\mu^e_t\}_{e \in \mathcal{E}_t}\}_{t \in T}$ must also be updated. Let the travel time of the edge $(v_o, v_d)$ be $\mu^{(v_o, v_d)}_t$. To generate the new distribution for the two edges, $(v_o, v_{\text{new}}^{(v_o, v_d)})$ and $(v_{\text{new}}^{(v_o, v_d)}, v_d)$, the value of $\mu^{(v_o, v_d)}_t$ is divided proportionally based on the relative lengths of the new edges. Let the $\gamma_t$ of spending time at $t$ through $v_o$ and $v_d$ as the measurement of the length ratio for $(v_o, v_{\text{new}}^{(v_o, v_d)})$ and $(v_{\text{new}}^{(v_o, v_d)}, v_d)$, respectively. The new values are computed as
\begin{align}
\mu^{(v_o, v_{\text{new}}^{(v_o, v_d)})}_{t+1} = \mu^{(v_o, v_d)}_t \cdot \gamma_t, \quad \mu^{(v_{\text{new}}^{(v_o, v_d)}, v_d)}_{t+1} = \mu^{(v_o, v_d)}_t \cdot (1-\gamma_t).
\label{eq:muupdate}
\end{align}
Travel time on the other edges remains the same.

So as well, let the original distribution of the edge $(v_o, v_d)$ be $\mathcal{C}^{(v_o, v_d)}_t$. we use a similar method to split complaint distributions
\begin{align}
\mathcal{C}^{(v_o, v_{\text{new}}^{(v_o, v_d)})}_{t+1} = \mathcal{C}^{(v_o, v_d)}_t \cdot \gamma_t, \quad \mathcal{C}^{(v_{\text{new}}^{(v_o, v_d)}, v_d)}_{t+1} = \mathcal{C}^{(v_o, v_d)}_t \cdot (1-\gamma_t).
\label{eq:Cupdate}
\end{align}
Complaint distribution on the other edges remains the same. Correspondingly, the output shape of the travel time predictor $\Xi$ and the complaints estimation $\Phi$.

In the meantime, the set of edge lengths also needs to be updated, including the deletion of $l^{(v_o,v_d)}$ and the appending of 
\begin{align}
l^{(v_o, v_{\text{new}})} = l^{(v_o, v_d)} \cdot \gamma_t, \quad l^{(v_{\text{new}}, v_d)} = l^{(v_o, v_d)} \cdot (1-\gamma_t).
\label{eq:lupdate}
\end{align}
Edge length on the other edges remains the same.

We aim to provide the patroller with greater flexibility when adapting during commuting. To achieve this, we split the neighboring edges of the current edge the patroller is on based on the complementary time spent on the total travel time. Define the neighbor edges of $v_{\mathbf{t}}$ as $\{(v_o=v_{\mathbf{t}},v_d=\overline{e})|\overline{e}\in \overline{N}(v_{\mathbf{t}}) = N(v_{\mathbf{t}})\setminus\{(v_{\mathbf{t}},a^*(W_{\mathbf{t}, \tau, |K|})),(a^*(W_{\mathbf{t}, \tau, |K|}),v_{\mathbf{t}})\}\}$, and the neighbor edges of $a^*(W_{\mathbf{t}, \tau, |K|})$ as $\{(v_o=a^*(W_{\mathbf{t}, \tau, |K|}),v_d=\underline{e})|\underline{e}\in \underline{N}(a^*(W_{\mathbf{t}, \tau, |K|}))= N(a^*(W_{\mathbf{t}, \tau, |K|}))\setminus\{(v_{\mathbf{t}},a^*(W_{\mathbf{t}, \tau, |K|})),(a^*(W_{\mathbf{t}, \tau, |K|}),v_{\mathbf{t}})\}\}$. The ratio of time spent is given by
\begin{align*}
    \gamma^{\overline{e}}_t = \frac{\tau-t +\mathbf{t}}{\mu_{\mathbf{t}}^{\overline{e}}}, \quad\text{and}\quad \gamma^{\underline{e}}_t = \frac{\tau - \mathbf{t}'+t}{\mu_{\mathbf{t}}^{\underline{e}}}.
\end{align*}
Following equations \eqref{eq:posnew}, \eqref{eq:Vupdate}, \eqref{eq:Eupdate}, \eqref{eq:muupdate}, \eqref{eq:Cupdate}, \eqref{eq:lupdate}, combined with, we get the updated $\mathcal{V}_{t+1}$,  $\mathcal{E}_{t+1}$:
\begin{align*}
    \mathcal{V}_{t+1}= \mathcal{V}_{t} \cup \{v^{(v_{\mathbf{t}}, a^*(W_{\mathbf{t}, \tau, |K|}))}_{\text{new}}\} \cup \{v^{\overline{e}}_{\text{new}}\}_{\overline{e}\in \overline{N}(v_{\mathbf{t}})} \cup \{v^{\underline{e}}_{\text{new}}\}_{\underline{N}(a^*(W_{\mathbf{t}, \tau, |K|}))},
\end{align*}

Let the sets of new edges be:
\begin{align*}
    \overline{\mathcal{E}}_{t+1}=\{(v_{\mathbf{t}}, v^{\overline{e}}_{\text{new}}), (v^{\overline{e}}_{\text{new}},\overline{e})\}_{\overline{e}\in \overline{N}(a^*(W_{\mathbf{t}, \tau, |K|}))}, \quad \underline{\mathcal{E}}_{t+1}=\{(v_{\mathbf{t}}, v^{\underline{e}}_{\text{new}}), (v^{\underline{e}}_{\text{new}},\overline{e})\}_{\underline{e}\in \underline{N}(a^*(W_{\mathbf{t}, \tau, |K|}))}.
\end{align*}
thus
\begin{align*}
\mathcal{E}_{t+1} = \{\mathcal{E}_t \setminus \{e\}_{e\in (v_{\mathbf{t}}, a^*(W_{\mathbf{t}, \tau, |K|})) \cup \overline{N}(v_{\mathbf{t}})\cup \underline{N}(a^*(W_{\mathbf{t}, \tau, |K|}))}\} \cup \overline{\mathcal{E}}_{t+1}\cup \underline{\mathcal{E}}_{t+1}.
\end{align*}

Therefore, the adapted patrol graph is represented as:
\begin{align*}
    \mathcal{G}_{t+1} = \langle \mathcal{V}_{t+1}, \mathcal{E}_{t+1} \rangle.
\end{align*} Then, we have updated the predictor of travel time and the complaint distributions $\Xi(W_{\mathbf{t}, \tau, |K|}), \Phi(W_{\mathbf{t}, \tau, |K|}, \{\hat{\mathcal{C}}_{e}\}_{e \in \mathcal{E}_{t+1}})$, respectively. Eventually, starting from the current time $ \mathbf{t}_{\text{new}} = t+1 $, the patroller solves the new patrolling problem $\left\langle W_{\mathbf{t}_{\text{new}}, \tau, |K|}, \{\hat{\mathcal{C}}^{e}_{\mathbf{t}_{\text{new}}}\}_{e \in \mathcal{E}_{\mathbf{t}_{\text{new}}}}, \Xi(\cdot), \Phi(\cdot), \rho = \{v_0 = v_{\mathbf{t}_{\text{new}}}, \mathcal{V} = \mathcal{V}_{\mathbf{t}_{\text{new}}}, \mathcal{E}= \mathcal{E}_{\mathbf{t}_{\text{new}}}, \mathcal{L}=\mathcal{L}_{\mathbf{t}_{\text{new}}}\} \right\rangle.$

\begin{algorithm}[H]
\caption{Online Adaptation Procedure}
\label{alg:online-adapt}
\begin{algorithmic}[1]
\Require DKW Shifting threshold $q$, Initial $\mathcal{G}_{\mathbf{t}}$, $\mathcal{L}_{\mathbf{t}}$, Prior $\{\hat{\mathcal{C}}^{e}_{\mathbf{t}}\}_{e}$, Predictor $\Xi$, Complaint $\Phi$, $\mathbf{t}$, $\mathbf{t}'$, $a^*(W_{\mathbf{t}, \tau, |K|})$

\For{$t = \mathbf{t}, \ldots, \mathbf{t}'$}
    \State Observe $o_t$ and update predictor as in \cite{ttcx}.
    \State Update $\{\hat{\mathcal{C}}_{e}(t)\}_{e}$ using Bayesian update \eqref{eq:update}.
    \State Compute $D_t^q$ by \eqref{eq:DKW_index}.
    \If{$D_t^q \neq 0$}
        \If{$t = \mathbf{t}'$} 
            \Comment{Inspecting}
            \State Solve the problem using Algorithm~\ref{alg:patrol-decision} starting from $t+1$.
        \ElsIf{Patroller on a node}
            \Comment{Inspecting}
            \State \textbf{break from for-loop}
            \State Solve the problem using Algorithm~\ref{alg:patrol-decision} starting from $t+1$.
        \Else
            \Comment{Commuting}
            \State Following the instruction in \Cref{sec:commuting_update}.
            \State \textbf{break from for-loop}
            \State Solve the problem using Algorithm~\ref{alg:patrol-decision} starting from $t+1$.
        \EndIf
    \EndIf
\EndFor
\end{algorithmic}
\end{algorithm}

Referring \eqref{eq:global_opt}, suppose we collect the trajectory of actions by TAMPA. Let $a_{\hat{\eta}_0}=v_0$ since it is the initial node of the patroller. From $t=0$ to $t=|T|$ denoted as $\hat{\mathcal{T}}=\{a_{\hat{\eta}_0}, a_{\hat{\eta}_1},\cdots, a_{\hat{\eta}_{\text{end}}}\}$. Thus, the global cost function of $\hat{\mathcal{T}}$ is denoted by

\begin{align}
 Q(\hat{\mathcal{T}}) =
  \sum_{g=0}^{\text{end}-1} 
  r_g\bigl(a_{\hat{\eta}_g}, a_{\hat{\eta}_{g+1}}\bigr).
  \label{eq:global_cost}
\end{align}

as the performance measurement of this action trajectory.

\section{Experimental Design}
\begin{figure}[h]
    \centering
    \includegraphics[width=0.8\linewidth]{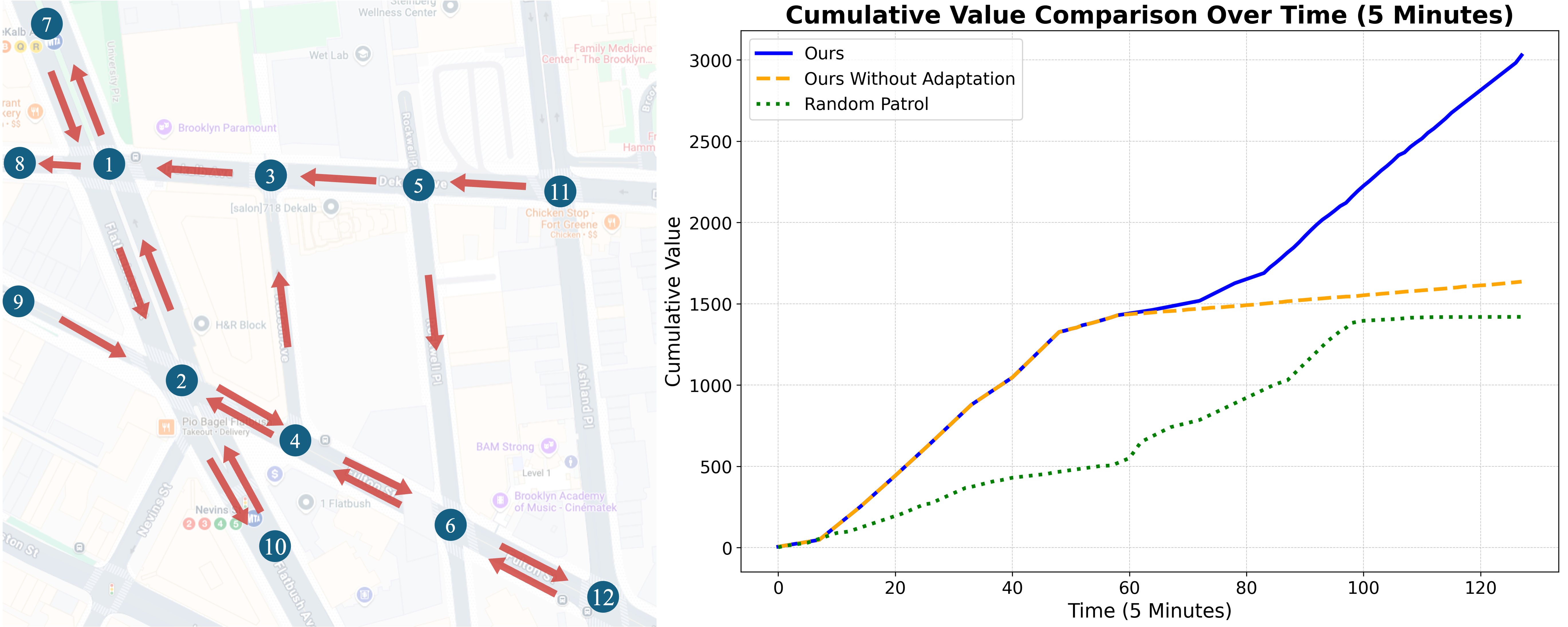}
    \caption{On the left, we illustrate the urban traffic network covered by the patroller, consisting of 12-node, 19-edge. On the right, we compare the cost in \eqref{eq:global_cost} of three patrolling strategies. The solid blue line represents our proposed method. The dotted yellow line represents solving \Cref{alg:patrol-decision} under a stationary setting (i.e., fixed graph). The dotted green line illustrates a baseline strategy in which the patrol vehicle selects its next destination node at random. Before 360 minutes (corresponding to the 75th time step under a 5-minute interval), the blue and yellow lines coincide because the complaint distribution remains unchanged. After the distribution shifts at 360 minutes, our method (blue line), which adapts to the new conditions, outperforms both the stationary approach and the random baseline. Overall, our method achieves an approximate $87.5\%$ improvement over the non-adaptive version of our approach and a $114.2\%$ improvement compared to a random patrol strategy.}
    \label{fig:patrol_exp_1}
\end{figure}

We apply TAMPA in an urban traffic network shown on the left side of \Cref{fig:patrol_exp_1} in Flatbush Avenue of Brooklyn, New York City. This location serves as a critical urban corridor with high traffic demand, particularly during mega-events hosted at venues such as the Barclays Stadium, home of NY NETS, which generates significant surges in vehicular and pedestrian flows and many concert events. The avenue’s complex traffic dynamics and high variability present an ideal environment for evaluating the effectiveness of the TAMPA in addressing real-time congestion and optimizing complaints management. Its strategic location and frequent disruptions due to large-scale events make it a representative case study for urban traffic systems requiring adaptive, scalable solutions. 

The dataset pivotal to training and evaluating TAMPA is generated using the Simulation of Urban Mobility (SUMO) \cite{sumo}, a versatile traffic simulation tool. The simulation models the Flatbush Avenue corridor, employing real-world data for calibration, including travel times, volumes, speeds, and surrogate safety metrics, to ensure realistic traffic dynamics \cite{sha2023calibrating}. Spanning 200 simulated days with one-second granularity, the dataset provides detailed minute-by-minute data on average volume, speed, and travel times for each road segment, enabling a comprehensive analysis of urban traffic patterns.

\begin{figure}[h]
    \centering
    \includegraphics[width=0.8\linewidth]{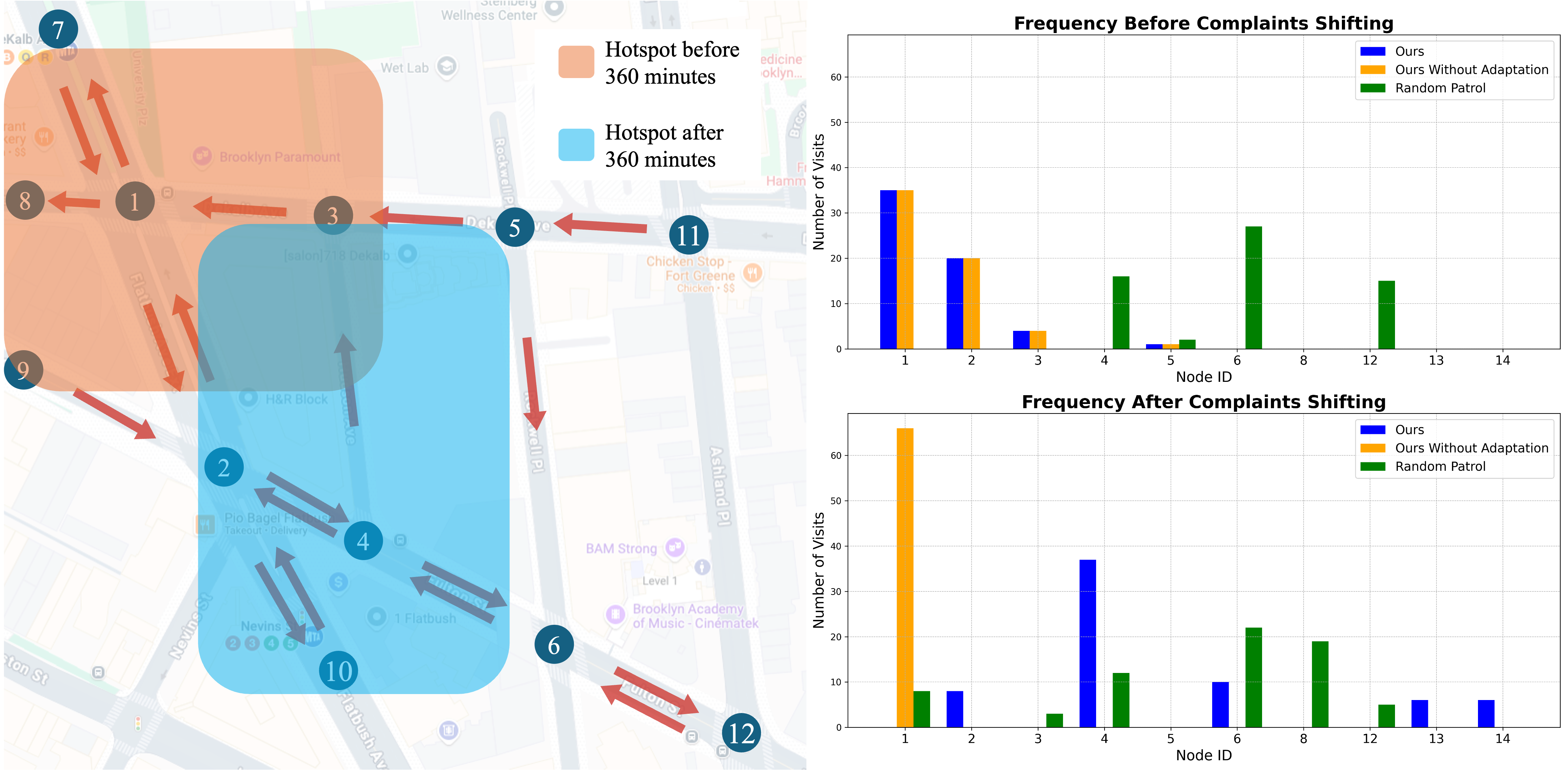}
    \caption{The left side illustrates how the complaint hotspot shifts from node $1$ to node $4$. On the right, we compare the visit frequencies of three different strategies before and after this shift. Notable, $13, 14$ are two newly added nodes that split edges $2-1$ and $2-4$. Our method adapts to the changing hotspot (spend almost $90\%$ time patrolling around node $4$ after the shifting), whereas the other two strategies fail to detect or respond effectively.}
    \label{fig:patrol_exp_2}
\end{figure}

To evaluate the performance of TAMPA, we conducted an initial experiment using a short simulation period with a time horizon of $|T|=700$ minutes. The initial patrol graph is a 12-node, 19-edge graph, as shown on the left of \Cref{fig:patrol_exp_1}. Accordingly, the state space consists of the $12$ nodes, with the initial state being node $5$. The action space at each node corresponds to the topology of the patrol graph. A special event was designed to occur 360 minutes into the simulation (randomly chosen to start at 360 minutes, though it could occur at any time).  As shown in \Cref{fig:patrol_exp_2}, before the special event, the distribution of complaints has a higher mean value around node $1$ (including its neighboring edges), designating node $1$ as the hotspot. After the special event, however, the mean value of the distribution around node $1$ decreases while it increases around node $4$, indicating that the hotspot shifts from node $1$ to node $4$.

We decide the fixed time interval $\tau$ as 8 minutes, referred the usual average travel time around the edges in the graph. By using the same predictor structure in \citep{ttcx}, and the 
complaints at each edge are computed as follows:
\begin{align*}
    \text{Noise: } & \quad n \sim \mathcal{N}(\mu=0.5, \sigma=0.2), \\
    \text{Edge Complaints: } & \quad c^e = \min\left(\max\left(\lfloor \mathcal{F}^e\cdot U(0,1) (1 + n) \rceil, 0\right), 30\right),
\end{align*}
where $n$ is the human behavior noise, $\mathcal{F}^e$ is the weight parameters for each edge, $U(0,1)$ represents a uniform distribution in the range $[0, 1]$, and $\lfloor \cdot \rceil$ denotes rounding to the nearest integer. 

Two comparison experiments illustrate two key improvements introduced by our design in \Cref{fig:patrol_exp_1}. First, an ablation study that removes the adaptation process on the traffic graph demonstrates that failing to respond to shifts in the complaints distribution can severely degrade patrol performance. Second, benchmarking against a random strategy, which selects the next destination at random, further underscores the benefits of our approach.

In \Cref{fig:patrol_exp_2}, we illustrate how the distribution of complaints shifts and compare the corresponding actions taken by three different strategies. When the hotspot of complaints moves from node 1 to node 4, neither the design lacking graph adaptation nor the random strategy can effectively detect and respond to these changes, underscoring the importance of an adaptive approach.

In summary, extensive simulation results illustrate that both removing the adaptation process and relying on random patrolling lead to substantially lower rewards, confirming the advantage of our adaptive approach in responding to dynamic traffic conditions in the presence of random accidents/incidents.

\section{Summary \& Future Works}
In this paper, we propose the Traffic Adaptive Moving-window Patrolling Algorithm (TAMPA), designed to mitigate the number of complaints during the management of high-impact events such as the FIFA World Cup. TAMPA integrates a traffic pattern predictor and a complaints estimator to balance immediate complaint coverage and routing constraints, using dynamic programming with moving horizons to derive optimal strategies adaptive to real-time observations. Since TAMPA utilizes empirical distributions to initialize the planning, which may significantly deviate from the observations. We design an adaptive method for adapting the patrol graph, based on the Dvoretzky–Kiefer–Wolfowitz testing. We prove that the suboptimality gap of the proposed TAMPA with respect to the optimal patrol plan in hindsight is upper bounded by the distribution shifts measured in total variation multiplied by the planning horizon. We validate our approach through real-world experiments in dynamic environments (around $87.5\%$ improvement compared with stationary strategy and $114.2\%$ improvement compared with random patrol strategy). 

Future efforts will focus on enhancing the proposed approach and exploring new methodologies to address current limitations. For instance, the flexibility of traffic network adaptation can be improved by incorporating automated mechanisms to remove less active nodes and edges. Moreover, integrating digital twin technology for both traffic pattern prediction and complaint estimation could introduce new capabilities and improve scalability. As part of these efforts, we are also working on developing a simulation model for the network surrounding MetLife Stadium, which will host the FIFA World Cup Final in 2026. This high-impact event provides a valuable opportunity to validate the approach in a complex and dynamic scenario. Lastly, experimental validation in diverse zones and scenarios, such as additional complaint hotspots and their dynamic movements, will remain a key focus moving forward.

\bibliographystyle{apalike} 
\bibliography{reference}
\end{document}